\documentclass{amsart}
\usepackage{amssymb,latexsym,graphics}
\usepackage[mathscr]{eucal}

\newtheorem{thm}{Theorem}[section]

\newtheorem{lem}[thm]{Lemma}

\theoremstyle{remark}
\newtheorem{rem}[thm]{Remark}

\theoremstyle{definition}

\newcommand{\lp}[2]{\Vert \, #1 \, \Vert_{#2}}
\newcommand{\td}[1]{\widetilde{#1}}
\newcommand{\bu}{\underline{u}}
\newcommand{\bL}{\underline{L}}
\newcommand{\RR}{\mathbb{R}}
\newcommand{\snabla}{ {\slash\!\!\!\! \nabla} }

\newcommand{\sgn}{\hbox{sgn}}

\newcommand{\ret}{\vspace{.3cm}}

\begin{document}

\title[Decay for Scalar Fields on Schwarzchild Space]
{Uniform Decay of Local Energy and the Semi-Linear Wave Equation 
on Schwarzchild
Space}
\author{Pieter Blue}
\address{Department of Mathematics
University of Toronto
100 St. George Street
Toronto, Ontario, M5S 3G3
Canada}
\email{pblue@math.toronto.edu}
\author{Jacob Sterbenz}
\address{Department of Mathematics,
University of California, San Diego (UCSD)
9500 Gilman Drive, Dept 0112
La Jolla, CA 92093-0112 USA}
\email{jsterben@math.ucsd.edu}
\thanks{
The second author would like thank MSRI and Princeton University, 
where a portion
of this research was conducted during the Fall of 2005.}
\subjclass{}
\keywords{}
\date{}
\commby{}


\begin{abstract}
We provide a uniform decay estimate for the local energy of general
solutions to the inhomogeneous wave equation on a Schwarzchild
background. Our estimate implies that such solutions have asymptotic
behavior $|\phi| =
O\left(r^{-1}(\big|t-|r^*|\big|^{-\frac{1}{2}}\right)$ as long as
the source term is bounded in the norm $(1-\frac{2M}{r})^{-1}\cdot(1
+ t + |r^*|)^{-1}L^1\big(H^3_\Omega(r^2dr^* d\omega)\big)$. 
In particular this gives
scattering at small amplitudes 
for non-linear scalar fields of the form $\Box_g\phi =
\lambda |\phi|^p\phi$ for all $2 < p$.
\end{abstract}

\maketitle

\section{Introduction}
In this paper our goal is to give a somewhat elementary discussion
of the global decay properties of general solutions to the  scalar
wave equation on the exterior of a Schwarzchild black hole. That is,
we consider the manifold with boundary:
\begin{equation}
        \mathcal{M} \ = \ \RR\times \,
    [2M,\infty) \times \mathbb{S}^2 \ , \label{S_space}
\end{equation}
with Lorentzian metric:
\begin{equation}
        ds \ = \ -(1 - \frac{2M}{r})\, dt^2 +
    (1-\frac{2M}{r})^{-1}dr^2 + r^2d\omega^2 \ , \label{S_metric}
\end{equation}
and we look at smooth functions $\phi$ which do not touch the
boundary of \eqref{S_space} for each \emph{fixed} value of the
parameter $t$ and which satisfy the inhomogeneous wave equation:
\begin{equation}
        \Box_g \phi \ = \
        \nabla^\alpha \partial_\alpha \phi \ = \ G . \notag
\end{equation}
The main question we would like to answer here is: How quickly
does the local
energy of the wave $\phi$ dissipate over compact sets in the $r$
coordinate, and how precisely does the dissipation depend on the
source $G$?\\

In the special case of Minkowski space, $M=0$, a quite satisfactory
answer to this question is known. Here one has the classical uniform
local energy decay estimate of C. Morawetz:
\begin{multline}
        \int_{\mathbb{R}^3\times \{t\}}\ \left( (1 + \bu^2)|L\phi|^2
    + (1 + u^2)|\bL \phi|^2 + (1 + u^2 + \bu^2)|\snabla\phi|^2
    + \frac{1 + \bu^2 + u^2}{r^2}|\phi|^2
    \right)\ dx \\
    \lesssim \ \left(\int_0^t\
    \lp{(1 + |\bu| + |u|)G(s)}{L^2(dx)}\ ds
    \right)^2 \ + \ \int_{\RR^3\times\{0\}}\ (1 +
    r^2)|\nabla_{t,x}\phi|^2 \ dx \
    . \label{class_mor}
\end{multline}
Here one sets $u=(t-r)$, $\bu=(t+r)$, $L=2\partial_{\bu}$, and
$\bL = 2\partial_u$.
For the original proof see the paper of Morawetz \cite{M_class}, and for an
alternative proof as well as many generalizations see the recent work
\cite{LS_MKG}.\\

The beauty of the estimate \eqref{class_mor} is that it gives one a
huge amount  of information about the global dispersive properties
of the function $\phi$. For one, it produces a \emph{pointwise} in
time decay of the $L^2_{loc}$ norm as well as the local energy.
What's more, this local decay is given in such a way that it is
clear what is happening on the whole of each time slice $t=const.$,
even very far away from the origin $r=0$. In fact, using Sobolev
embeddings and \emph{only} rotations, the Morawetz estimate is good
enough to provide uniform decay at the rate of $(1+t)^{-1}$.
However, perhaps the most important property of the estimate
\eqref{class_mor} is that it turns out to be incredibly useful when
dealing with non-linear problems. This is because it places very
simple conditions on the source term $G$, the kind which are
relatively straight forward to recover in bootstrapping arguments
given the form of the left hand side of \eqref{class_mor}. For
instance, \eqref{class_mor} makes dealing with the global existence
problem for small amplitude non-linear scalar fields of the form:
\begin{equation}
        \Box \phi \ = \ \lambda \, |\phi|^{p}\phi \ , \label{NLW}
\end{equation}
essentially trivial in the case where $2 < p$. All one has to do is
to combine   \eqref{class_mor} with appropriately localized Sobolev
embeddings to yield the decay estimate:
\begin{equation}
        |\phi| \ \lesssim \ (1 + r)^{-1}\cdot\min\{\frac{r^\frac{1}{2}}
        {1 + |t-r|} \ , \ \frac{1}{1 + |t-r|^\frac{1}{2}}\}
         \ , \notag
\end{equation}
which is enough to integrate the nonlinearity on the right hand side
of \eqref{NLW} when it appears on the right hand side of
\eqref{class_mor}. In fact, if one takes into account characteristic
estimates of the form \eqref{class_mor}, see again \cite{LS_MKG},
then it is possible to push the exponent $p$ to certain values
$p\leqslant 2$. We will not discuss such refinements here.\\

What we will show here is that for the more general case of $M\neq
0$, an estimate which is essentially of the form \eqref{class_mor}
holds in the case of general Schwartz (on each fixed time slice)
functions $\phi$. The proof we give is a relatively straight forward
integration by parts, similar in spirit  to how one proves
\eqref{class_mor}. In the final section of the paper, we indicate
how our estimates can be used to give a short proof of global
existence and decay for non-linear scalar fields of
the form \eqref{NLW} when $2 < p$.\\

To state our main theorem, we will use the following notation. We
first reparametrize the radial variable in the usual way:
\begin{equation}
        r^* \ = \ r + 2M\, \ln(r-2M) \ , \label{RW_r*}
\end{equation}
and then introduce the optical functions and null-generators for the
coordinates $(t,r^*,\omega)$:
\begin{subequations}\label{defs}
\begin{align}
        \bu \ &= \ (t+r^*) \ ,
        &u \ &= \ (t-r^*) \ , \label{defs1} \\
        L \ &= \ \partial_t + \partial_{r^*} \ ,
        &\bL \ &= \ \partial_t - \partial_{r^*} \ . \label{defs2}
\end{align}
\end{subequations}
We will prove that:\\

\begin{thm}[Uniform Local Energy Decay for the
Scalar Wave Equation on Schwarzchild Space]\label{main_thm} Let
$(t,r^*,\omega)$ be the coordinates (as defined above) on the
manifold $\mathcal{M}$ \eqref{S_space} with metric \eqref{S_metric}.
Let $\phi$ be a smooth function compactly supported on each
hypersurface $t=const.$ and set:
\begin{equation}
        \Box_g \phi \ = \ G \ . \notag
\end{equation}
Then one has the following global estimate:
\begin{multline}
        \int_{\RR\times\mathbb{S}^2\times\{t\}}\ \Big((1 + \bu^2)(L (r\phi))^2
    + (1 + u^2)(\bL (r\phi))^2 \\
    + (1 + \bu^2 + u^2)(1 - \frac{2M}{r})\cdot\left[
    \frac{1}{r^2}|\snabla_\omega(r\phi)|^2 +
    \frac{M}{r^3}(r\phi)^2\right]
    \Big)\ dr^* d\omega \\ \lesssim \
    \left(\int_0^t \ \lp{(1 + |\bu| + |u|)
    (1 - \frac{2M}{r} )\, r\cdot
    (\sqrt{1-\Delta_{sph}}G)(s)}{L^2(dr^*d\omega)}\ ds
    \right)^2 \\
    + \int_{\RR\times\mathbb{S}^2\times\{0\}}\
    (1 + (r^*)^2) \ \Big[\big|\nabla_{t,r^*}\sqrt{1-\Delta_{sph}}(r\phi)\big|^2
     \\
     + \ (1 - \frac{2M}{r}) \left(
    \frac{1}{r^2}\big|\snabla_\omega(\sqrt{1-\Delta_{sph}}\, (r\phi))\big|^2 +
    \frac{M}{ r^3}(\sqrt{1-\Delta_{sph}}(r\phi))^2
     \right) \Big]\ dr^*d\omega
    \ , \label{morawetz}
\end{multline}
where the implicit constant  depends only on the mass $M$. 
Here $\Delta_{sph}$ is
the Laplacean in the angular variable $\omega$, and $\snabla_\omega$
is the associated gradient.\footnote{This should not
be confused with $\snabla$ from
line \eqref{class_mor} which is the covariant derivative on spheres of
radius $r$. Of course these two only differ by the factor of $r^{-1}$.}
\end{thm}\ret

\begin{rem}
Let us first give a heuristic summary of the content of the estimate
\eqref{morawetz} and how it contrasts to the situation of flat space
$M=0$. Roughly speaking, the Schwarzchild space can be split into
three pieces where one sees qualitatively  different behavior in
solutions to the wave equation.\\

The first region is very close to the boundary $r=2M$. For the static
space--time we are considering, this is quite easy to 
understand. Here wave propagation looks
essentially trivial in that one has  $\phi\sim F_1(t+r^*,\omega)$
for some smooth decaying function $F_1$ on the space $\RR\times
\mathbb{S}^2$. That is, wave propagation near the boundary $r=2M$ is
essentially just transport in the variable $\bu = t+r^*$. The caveat
is that this variable is the only one in which it is possible to get
decay for this region because $F_1(0,\omega)$ does not need to be
small. \\

We would like to call the readers attention to the fact
that the precise decay of the function
$F_1$ in the first variable seems to be quite
delicate issue, and we will only
obtain $|F_1(\bu,\omega)| \lesssim (1+|\bu|)^{-\frac{1}{2}}$. Of
course, this is all one should expect given that the right hand side
of \eqref{morawetz} is consistent with initial data decaying at this
rate. Since this estimate does not ask for a lot of information,
which is actually its strength in treating non-linear problems, it
does not give a lot of information in return. For a much more
precise asymptotic in the case of spherical symmetry, and for the
more difficult case of dynamic space-times (for the spherically
symmetric coupled scalar field), we refer the reader to the very 
deep recent work of Dafermos-Rodnianski (see \cite{DR_1}--\cite{DR_2}).\\

The second region is in the ``far exterior'' $t\lesssim r^*$ where
things look essentially flat. This is also fairly easy to
understand. Here one expects that things look very similar to the
well known case of Minkowski space.\\

The third region is ``the boarder'' close to $r=3M$, which in
Regge--Wheeler coordinates \eqref{RW_r*} we extend to 
the region $|r^*|\leqslant
\frac{1}{2} t$. This is by far the most difficult region to
understand, and where one looses regularity in the estimate
\eqref{morawetz}. This loss of regularity is in sharp contrast to
the estimate \eqref{class_mor} in the case of Minkowski space, and
is also something one sees only by looking at the non-spherically
symmetric (functions) situation. What is happening here is that the
geometry is pulling the radiation apart into the two regions just
mentioned, and there is a danger that this ``splitting'' could allow
some fairly large residual portion of the radiation to linger for a
long amount of time in the transition region $r^*\sim 0$.\\

Now it turns out that this effect can only happen (and it does
happen!) if the wave $\phi$ has a very high angular momentum. In
this case it can concentrate on a very small set for each fixed time
in the $\omega$ variable, and it will essentially rotate around the
sphere $r=3M$ for a long while before dispersing. This behavior 
can be understood by observing 
that null geodesics tangent to the surface $r=3M$ will remain tangent 
to this surface \cite{MTW}, and that high angular momentum solutions 
to the wave equation will closely follow these geodesics for a long 
period of time before dispersing.\\

This slow dispersion can also be understood by conformally 
changing the metric \eqref{S_metric} by the factor $(1-\frac{2M}{r})$. 
Since the coefficient of $dt^2$ is constant on this 
new manifold, the corresponding wave equation describes the 
time evolution of a wave on a three dimensional Riemannian
manifold with 
metric given by the spatial portion of the confromal Lorentzian metric. 
A simple calculation shows that this Riemannian manifold
has a totally geodesic 
sphere (and hence closed geodesics) at the value $r=3M$.
Now, the original wave equation is equivalent to the wave
equation with respect to the
conformal metric modulo a smooth potential. For very
high frequencies this potential cannot compete with the principle
part of the conformal wave operator, so it is not difficult to
construct coherent state solutions which concentrate near the closed
geodesics sitting at $r=3M$ for a long amount of time. Therefore,
from this point of view, one should look at \eqref{morawetz} as a
sort of ``cheap'' dispersive estimate, and it is well known that
such estimates loose regularity when the underlying geodesic flow is
not well behaved in  the sense of spreading of the classical
trajectories.\\

We further remark here that the nature of the surface and geodesics at 
$r=3M$ can be a 
little confusing to discuss in the relativistic terminology. The 
null geodesics at $r=3M$ which orbit the black hole form a helix 
in four dimensional space-time with an axis in the time direction. 
Although their projections onto the three dimensional coordinates 
$(r,\omega)$ is closed, because the $t$ coordinate is constantly 
increasing they are not closed null geodesics. Also, while it is true
that the hypersurface $r=3M$ is foliated by null geodesics 
it is not itself a null hypersurface, which is one with a null normal 
direction (see \cite{HawkingEllis}). Note that a normal to $r=3M$ is 
$\partial_{r^*}$ which is space like.\\

As far as our analysis is concerned,
the presence of null geodesics at $r=3M$ manifests itself 
through trapping terms which are positive for $r\sim 3M$. 
For a wave equation with a potential $Q$, we refer 
to $\vec{x} \cdot \nabla Q + 2Q$ as the trapping term. This 
expression appears as a contribution governing the growth of the 
conformal (Morawetz) energy. It can be seen as the main ``error''
which is generated by the divergence of the conformal energy density,
and is given by the first two terms on the right hand side of equation
\eqref{K0_div_iden} below. This identity was first derived   
using a some what different formalism in the dissertation 
of the first author (see \cite{PB_T}). 
\end{rem}\ret

\begin{rem}
Our proof of the bound \eqref{morawetz} will be very general in the sense
that we derive it from a fairly generic family of estimates that
holds for 1-D wave equations with ``strongly repulsive'' potentials.
It is to be hoped that this procedure can be used to accommodate
other situations, such as higher spin equations on Schwarzchild
space and possibly other space-times where spherical harmonic
decompositions make sense. We will leave these discussions to
further work. \\

The approach we take here is based on the previous works 
\cite{SB_1}--\cite{SB_2} which proved space-time Morawetz type estimates  
on Schwarzschild-like manifolds, and the thesis \cite{PB_T}
which proved a version of the conformal (Morawetz) energy 
estimate \eqref{morawetz} with growing right hand side. 
In the estimate contained in \cite{PB_T},  the trapping term 
(described previously) appears integrated in space-time 
against the quantity $t(\phi)^2$, where $\phi$ is the scalar field.  
If the factor of $t$ were not present and if the field $\phi$ were 
restricted to a single spherical harmonic, then the Morawetz 
estimate from \cite{SB_1}
would be sufficient to control this space-time integral.
However, due to the fact that the reduction to individual spherical 
harmonics leads to trapping terms which grow quadratically according to  
the angular frequency, both an additional angular derivative \emph{and}
the factor of $t$ must be controlled.\\ 

In this paper, we present a simple argument which allows one to absorb
the trapping term with the factor of $t$, and hence prove \eqref{morawetz}. 
In the dissertation \cite{PB_T} and the forthcoming 
work \cite{BS_3}, a more involved phase space analysis is used to reduce the 
loss of angular regularity in the space-time Morawetz estimate and in the 
analogue of \eqref{morawetz} to only
$\epsilon$ powers of the operator $1-\Delta_{sph}$.  
We leave the combination of these two techniques to future work.

\end{rem}\ret

\noindent
{\em{Acknowledgments:}} First and foremost we would each like to thank
our advisors Matei Machedon and Avi Soffer for introducing us to
the problem of decay estimates for the wave equation on 
Schwarzchild space, and for generously sharing their ideas with us. 
Without them this paper would, quite simply, not exist.
We would also like to thank Hans 
Lindblad for reading a preliminary version of this article and 
giving several useful comments, as well providing a good deal of
support and overall encouragement.
Special thanks also to 
Jason Metcalfe, Igor Rodnianski, and Chris Sogge
for many helpful conversations along the way.

\ret

\section{Preliminary Setup}
In this section we will set up some preliminary notation and ideas
from one dimensional wave equations on Minkowski space. This
material is for the most part entirely standard, and we make no
claim of originality for the basic concepts. Now, it turns out that
Theorem \eqref{main_thm} is actually a special case of a family of
estimates which holds in the following general situation. We
consider 1-D wave equations of the form:
\begin{equation}
        \Box \psi - Q(x)\psi \ = \ H \ , \label{1d_wave}
\end{equation}
where $\Box = -\partial_t^2 + \partial_x^2$ and $Q(x)$ is some smooth real
valued function which we assume is general for the time being. When
the source term $H$ vanishes the field \eqref{1d_wave} comes from a
Lagrangian with energy momentum tensor:
\begin{equation}
        T_{\alpha\beta}[\psi] \ = \
    \partial_\alpha\psi\partial_\beta\psi -
    \frac{1}{2}g_{\alpha\beta}
    \left(\partial^\gamma\psi\partial_\gamma\psi
    + Q(x)\cdot(\psi)^2\right) \ . \label{1d_em}
\end{equation}
A quick calculation using the equation \eqref{1d_wave} shows that one
has the divergence identity:
\begin{equation}
        \partial^\alpha T_{\alpha\beta}[\psi] \ = \
    -\frac{1}{2}\partial_\beta(Q)\cdot (\psi)^2
    + H\cdot\partial_\beta \psi \ . \label{gen_divergence}
\end{equation}
Also,  the trace identity:
\begin{equation}
        g^{\alpha\beta}T_{\alpha\beta}[\psi] \ = \
    -Q(x)\cdot(\psi)^2 \ . \notag
\end{equation}
follows immediately, where $g=diag(-1,1)$ is the 1-D Minkowski
metric.\\

The utility of the tensor \eqref{1d_em} is that is keeps track of how
the field \eqref{1d_wave} reacts to the flow of various vector-fields
$X=X^\alpha\partial_\alpha$ on $\RR\times \RR$. In general, we form
the momentum density associated to $X$:
\begin{equation}
        {}^{(X)}P_\alpha \ = \
    T_{\alpha\beta}X^\beta \ , \label{momentum_density}
\end{equation}
and we compute from \eqref{gen_divergence} the divergence:
\begin{equation}
        \partial^\alpha\, {}^{(X)}P_\alpha \ = \
    -\frac{1}{2}X(Q)\cdot(\psi)^2 + \frac{1}{2}
    T_{\alpha\beta}{}^{(X)}\pi^{\alpha\beta}
    + H\cdot X(\psi) \ , \label{div_iden}
\end{equation}
where ${}^{(X)}\pi$ is the deformation tensor:
\begin{equation}
        {}^{(X)}\pi_{\alpha\beta} \ = \
    \partial_\alpha X_\beta + \partial_\beta X_\alpha \ . \notag
\end{equation}
In the next section we will use this setup to prove the following
general 1-D uniform local energy decay estimate:\\

\begin{thm}[1-D Morawetz Estimate for Positive Strongly
Repulsive Potentials] \label{main_thm_red} Let $\psi$ be a function
on $(1+1)$ Minkowski space which is compactly supported for each
fixed value of the time variable $t$. Suppose $\psi$ satisfies the
equation \eqref{1d_wave} for some smooth real valued function $Q(x)$
which satisfies all of the following conditions:
\begin{align}
        0 \ &\leqslant \ Q \ , && &\hbox{(Positivity)}&
    \label{rep1}\\
        0 \ &\leqslant \ -\, x \partial_x(Q) \  && &\hbox{(Repulsive 1)}&
    \label{rep2}\\
    x \partial_x(Q)(x) + 2Q(x) \ &\leqslant \
    - C\sgn(x)\, \partial_x(Q)(x) \ ,
    &x\notin \mathcal{B}_1& \ ,  &\hbox{(Repulsive 2)}&
    \label{rep3}\\
    x \partial_x(Q)(x) + 2Q(x) \ &\leqslant \
     C |x|^{-1} Q(x) \ ,
    &x\notin \mathcal{B}_2& \ ,  &\hbox{(Homogeneity)}&
    \label{rep4}\\
    (1+\lambda^2)x^2 \ &\leqslant \ - Cx\partial_x(Q)(x) \ ,
    &x\in 2\mathcal{B}_1& \ ,  &\hbox{(Critical Point)}&
    \label{rep5}\\
    C^{-1} \ \leqslant \ Q(x) \ &\leqslant \
    C(1+\lambda^2) \ , &x\in 2 \mathcal{B}_1& \ ,
    &\hbox{(Local Bounds)}& \label{rep6}
\end{align}
where $C$ and $\lambda$ are fixed non-negative constants with
$C$ strictly positive, and the
$\mathcal{B}_i$ are compact intervals containing the origin. Then
the following uniform local energy decay estimate of Morawetz type
holds:
\begin{multline}
        \int_{\RR\times\{t\}}\
    \left((1+\bu^2)(L \psi)^2 +
    (1 + u^2)(\bL \psi)^2 + (1 + \bu^2 + u^2)Q(x)\cdot
    (\psi)^2 \right)\ dx \\ \\
    \lesssim \ (1+\lambda^2)\underline{E}(\psi(0))\\ + \
    \int_0^t \lp{(1 + |\bu| + |u|)(1+\lambda)H(s)}{L^2(dx)}
    \cdot\lp{(1+\lambda)(|\nabla_{t,x}\psi| + 
    Q^\frac{1}{2}\cdot|\psi| )(s)}{L^2(dx)}
    \ ds  \ , \\
    + \ \int_0^t \lp{(1 + |\bu| + |u|)\, H(s)}{L^2(dx)}
    \cdot\underline{E}^\frac{1}{2}(\psi(s))\ ds \ .
    \label{1D_mor}
\end{multline}
Here we have set:
\begin{align}
        \bu \ &= \ t + x \ ,
    &u \ &= \ t-x \ , \notag\\
    L \ &= \ \partial_t + \partial_x \ ,
    &\bL \ &= \ \partial_t - \partial_x \ , \notag
\end{align}
and:
\begin{equation}
        \underline{E}(\psi(s)) \ = \
    \int_{\RR\times\{s\}}\
    \left((1+\bu^2)(L \psi)^2 +
    (1 + u^2)(\bL \psi)^2 + (1 + \bu^2 + u^2)Q(x)\cdot
    (\psi)^2 \right)\ dx \ . \label{mor_energy}
\end{equation}
The implicit constant in estimate \eqref{1D_mor} depends only on the
constant $C$ and the size of the two intervals $\mathcal{B}_i$, and
\emph{not} on the value of $t$ or $\lambda$ or any other property of
the potential $Q(x)$ then those listed above.
\end{thm}\ret

\subsection{The case of Schwarzchild Space}
Before moving on to prove the estimate \eqref{1D_mor}, let us first
briefly indicate how this can be used to prove the bound
\eqref{morawetz}. In the $(t,r^*,\omega)$ coordinates one writes the
wave operator $|g|^{-\frac{1}{2}} \partial_\alpha g^{\alpha\beta}
|g|^\frac{1}{2}\partial_\beta$ as:
\begin{equation}
        (1-\frac{2M}{r})^{-1}\left(-\partial_t^2\phi + 
        r^{-2}\partial_{r^*} (r^2
        \partial_{r^*}\phi)\right) + \frac{1}{r^2}\Delta_{sph}\phi
        \ = \ G \ . \notag
\end{equation}
Here $\Delta_{sph}$ is the Laplacean in the $\omega$ variable.
Introducing now the quantities $\overline{\psi} = r\phi$ and
$\overline{H} = (1-\frac{2M}{r}) r G$ this last line becomes:
\begin{equation}
        -\partial_t^2 \overline{\psi} +
        \partial^2_{r^*}\overline{\psi} - r^{-1} \partial_{r^*}^2(r)
        \overline{\psi} + \frac{(1-\frac{2M}{r})}{r^2}
        \Delta_{sph}\overline{\psi} \ = \ \overline{H} \ \
        . \label{conf_eq}
\end{equation}
We now follow the usual procedure of projecting this equation onto
individual spherical harmonics. Since all of  our estimate are both
$L^2$ and rotation invariant, there is absolutely no harm in doing
this. We write:
\begin{equation}
        \overline{\psi} \ =\
        \sum_{\lambda,i}\ \psi_{\lambda,i} Y^i_\lambda \ , \notag
\end{equation}
where the $Y^i_\lambda$ form an orthonormal  basis for the space
$\Delta_{sph}Y_\lambda = -\lambda^2 Y_\lambda$. On each harmonic the
equation \eqref{conf_eq} becomes:
\begin{equation}
        -\partial_t^2 \psi_{\lambda,i} +
        \partial^2_{r^*}\psi_{\lambda,i} - r^{-1} \partial_{r^*}^2(r)
        \psi_{\lambda,i} - \frac{\lambda^2 (1-\frac{2M}{r})}{r^2}
        \psi_{\lambda,i} \ = \ H_{\lambda,i} \ \
        . \label{ang_conf_eq}
\end{equation}
Dropping the $(\lambda,i)$ indices, labeling $r^*=x$, and using the
notation:
\begin{equation}
        Q_\lambda \ = \ r^{-1} \partial_{r^*}^2(r)
        + \frac{\lambda^2 (1-\frac{2M}{r})}{r^2} \ , \label{l_pot}
\end{equation}
equation \eqref{ang_conf_eq} becomes:
\begin{equation}
        -\partial_t^2 \psi + \partial_x\psi - Q_\lambda(x)\psi \ =
        \ H \ . \notag
\end{equation}
We now wish to apply the estimate \eqref{1D_mor} to each of these
equations, after we apply a spatial translation by a quantity
$x_0(\lambda)$ which will be determined in a moment. The resulting
family of estimates can then be safely added to obtain the full
estimate \eqref{morawetz} as long as one can produce a \emph{single}
point $x_0(\infty)$ such that $x_0(\lambda)\to x_0(\infty)$, and a
\emph{single} set of objects $(C,\mathcal{B}_1,\mathcal{B}_2)$ such
that the conditions \eqref{rep1}--\eqref{rep6} hold for
$(C,\mathcal{B}_1+x_0(\lambda),\mathcal{B}_2+x_0(\lambda))$.
Luckily, for the family of potentials \eqref{l_pot} where $0
\leqslant \lambda$ is any real number, this is easy to show. The
reader should keep in mind that the reason this is possible is that
the conditions \eqref{rep1}--\eqref{rep6} are not really \emph{size}
conditions on the potential $Q_\lambda$, but are actually conditions
on the \emph{sign} of $Q_\lambda$ and its first derivative. This
type of condition is very stable under multiplication by large
positive constants, so it is not hard to get uniform behavior for
large $\lambda$. We will state this result as follows:\\

\begin{lem}\label{sch_lem}
Let $Q_\lambda$ be the family of potentials defined on line
\eqref{l_pot} above, and set $x_0(\infty)=3M$. Then there exists a
constant $C$, a pair of sets $\mathcal{B}_1,\mathcal{B}_2$, and a
family of points $x_0(\lambda)\to x_0(\infty)$ such that the
potentials $Q_\lambda\left(x+x_0(\lambda)\right)$ satisfy the
conditions \eqref{rep1}--\eqref{rep6} for the triple
$(C,\mathcal{B}_1+x_0(\lambda),\mathcal{B}_2+x_0(\lambda))$. All of
these objects are completely determined by the value of $M$.
\end{lem}\ret

\begin{proof}
First, notice that condition \eqref{rep1} is immediate. Next, recall
that in the current notation we have $x = r(x) + 2M\ln(r(x) - 2M)$.
We now write:
\begin{equation}
        Q_\lambda(x) \ = \ (1-\frac{2M}{r})\left(\frac{2M}{r^3}
        + \frac{\lambda^2}{r^2}\right) \ . \label{new_Qlambda_def}
\end{equation}
The proof centers around showing that $Q_\lambda$ has an isolated
critical point. We compute the first derivative:
\begin{align}
         Q_\lambda' \ &= \ \frac{2M}{r^2}(1-\frac{2M}{r})\left(\frac{2M}{r^3} +
         \frac{\lambda^2}{r^2}\right)\ - \  (1-\frac{2M}{r})^2\left(
        \frac{6M}{r^4} + \frac{2\lambda^2}{r^3} \right) \ , \notag\\
        &= \ -\frac{2}{r^5}(1-\frac{2M}{r})\cdot \left[\lambda^2 r^2 -
      3M(\lambda^2-1)r -
        8M^2\right]\ . \label{main_poly}
\end{align}
The polynomial on the right hand side of this last expression has
exactly one root for positive values of $r$. This is given by the
quadratic formula:
\begin{equation}
        r(\lambda) \ = \ \frac{3M(\lambda^2-1) +
        M\sqrt{9(\lambda^2-1)^2 + 32\lambda^2}}{2\lambda^2} \
        . \label{r_lambda_def}
\end{equation}
We now show that this positive root is trapped inside the interval
$[\frac{8M}{3},3M]$. Since it is clear from \eqref{r_lambda_def} that
asymptotically $r(\lambda)\to 3M$, it suffices to show that
$r(\lambda)$ is an increasing function for $0\leqslant \lambda$. This
follows at once from differentiating the polynomial on line
\eqref{main_poly} with respect to the parameter $\lambda$ which
yields the identity:
\begin{equation}
        \dot{r}(\lambda) \ = \ \frac{6M\lambda r - 2\lambda r^2}
    {2\lambda^2 r - 3M(\lambda^2-1)}  \ . \notag
\end{equation}
A simple calculation shows that this quantity is indeed positive whenever
$r\in[\frac{8}{3}M,3M]$.
Therefore, we shall pick our sequence of points $x_0(\lambda)$
according to the rule $r(x_0(\lambda))=r(\lambda)$. This immediately
gives the condition \eqref{rep2} for the family of translated
potentials $Q_\lambda(x + x_0(\lambda))$. Also, note that the
pointwise bound \eqref{rep6} is immediate for any compact interval
$\mathcal{B}_1$.\\

We now show the critical point behavior \eqref{rep5}. This boils
down to the fact that the polynomial on line \eqref{main_poly} has a
simple root at $r(\lambda)$. In fact, taking the second derivative
of the potential $Q_\lambda$ with respect to $x$ and evaluating at
the point $x_0(\lambda)$ we have that:
\begin{equation}
        Q''_\lambda(x_0(\lambda) )
    \ = \ -\frac{2}{r^5(\lambda)}
    (1-\frac{2M}{r(\lambda)})^2 \cdot \left[2\lambda^2 r(\lambda) -
      3M(\lambda^2-1) \right] \ . \notag
\end{equation}
Notice that this quantity never vanishes, and is $O(\lambda^2)$ as
$\lambda\to\infty$, so one has \eqref{rep5} for \emph{any} compact
set $\mathcal{B}_1$ given a suitable constant $C$, independent of the
value of $\lambda$.\\

It remains for us to show the ``strongly repulsive'' conditions
\eqref{rep3}--\eqref{rep4} hold for a uniform constant $C$ and pair of
sets $\mathcal{B}_i$. This follows from direct inspection of the
formulas \eqref{new_Qlambda_def} and \eqref{main_poly}. We consider
the cases of $x\to\pm\infty$ separately. In the case of $x\to\infty$
we also have that $r\to\infty$, and we have the two asymptotic
formulas (with uniform constants in $\lambda$ depending only on the
mass $M$):
\begin{align}
        (x-x_0(\lambda))\cdot Q'_\lambda(x) \ &= \
        -\frac{2\lambda^2}{r^2} - \frac{6M}{r^3} +
        O(\frac{\lambda^2}{r^3}) + O(\frac{1}{r^4})
    + \{\hbox{something negative}\}\ , \notag\\
    2Q_\lambda(x) \ &= \ \frac{4M}{r^3} + \frac{2\lambda^2}{r^2}
    + \{\hbox{something negative}\}\ . \notag
\end{align}
Notice that the \{something negative\} terms on the right hand side of
the first line above contain logarithmic
corrections of the form $\ln(x)/x^4$ and $\lambda^2 \ln(x)/x^3$,
which come from the second summand on the 
right hand side of \eqref{RW_r*}. It is important these come with 
a good sign.
Now, combining these last two lines we have that as $x\to \infty$:
\begin{equation}
        (x-x_0(\lambda))\cdot Q'_\lambda(x) + 2Q_\lambda(x)
        \ \leqslant  \  - \frac{2M}{r^3}  +
        O(\frac{\lambda^2}{r^3}) + O(\frac{1}{r^4}) \ . \notag
\end{equation}
This is enough to imply \eqref{rep3}--\eqref{rep4} because as
$x\to\infty$ we also have the following strict lower bounds:
\begin{align}
        \frac{\frac{1}{2}\lambda^2}{x^3}
    \ &\leqslant \ -Q' \ ,
    &\frac{\frac{1}{2}\lambda^2}{x^3}
    \ &\leqslant \ \frac{1}{x} Q_\lambda(x)
    \ . \notag
\end{align}\ret

Finally, we deal with the bounds \eqref{rep3}--\eqref{rep4} as $x\to
-\infty$. This is even easier to treat. Notice simply that both $Q$
and $\partial_x(Q)$ are
$O\left((1+\lambda^2)(1-\frac{2M}{r})\right)$, while the factor
$(x-x_0(\lambda))$ goes to $-\infty$. This means that the first term
on the left hand side of both  \eqref{rep3}--\eqref{rep4} is a very
large negative multiple of the second. Therefore, the bounds
\eqref{rep3}--\eqref{rep4} are trivial because the left hand side is
asymptotically negative. This completes our demonstration of Lemma
\ref{sch_lem}.
\end{proof}\ret

To wrap things up for this section, let us just mention briefly how
one can sum the estimate \eqref{1D_mor} over the angular frequency
localized pieces $\psi_{\lambda,i}$. The key thing here is that the
estimate \eqref{1D_mor} has been set up in such a way that one can
use the Cauchy-Schwartz inequality to deal with the terms on 
the right hand side
of \eqref{1D_mor}.  Specifically, summing this bound over
$(\lambda,i)$ indices and using the fact that the $\{Y^i_\lambda\}$
form an orthonormal system on the sphere $\mathbb{S}^2 $ we have
that:
\begin{multline}
        \sup_{0\leqslant s \leqslant t}\ \underline{E}(\phi(s)) \ = \
     \sup_{0\leqslant s \leqslant t}\
     \int_{\RR\times\mathbb{S}^2 \times\{s\}}\
    \Big((1+\bu^2)(L (r\phi) )^2 +
    (1 + u^2)(\bL (r\phi) )^2\\
    + (1 + \bu^2 + u^2)\cdot\big(\frac{|\snabla_\omega (r\phi)|^2}{r^2}
     + \frac{M}{r^3}
    (r\phi )^2\big) \Big)\ dr^* d\omega \\ \\
    \lesssim \ \sum_{\lambda,i}\ \Big[\ (1+\lambda^2)
    \underline{\underline{E}}(0)(\psi_{\lambda,i})\\ +
    \int_0^t \lp{(1 + |\bu| + |u|)(1+\lambda)
      \overline{H}_{\lambda,i}(s)}{L^2(dr^*)}
    \cdot\lp{(1+\lambda)(|\nabla_{t,r^*}\psi_{\lambda,i}| +
    (1-\frac{2M}{r})^\frac{1}{2}(\frac{\lambda^2}{r^2} +
    \frac{M}{r^3})^\frac{1}{2}
    |\psi_{\lambda,i}|
    )(s)}{L^2(dr^*)}
    \ ds  \  \\
    + \int_0^t \lp{(1 + |\bu| + |u|)\, 
      \overline{H}_{\lambda,i}(s)}{L^2(dr^*)}
    \cdot\underline{\underline{E}}^\frac{1}{2}
    (\psi_{\lambda,i}(s))\ ds \ \Big] \
    , \label{big_guy}
\end{multline}
where we are defining:
\begin{align}
        \underline{\underline{E}}(\psi_{\lambda,i}(s)) \ = \
        &\int_{\RR\times\{s\}}\ \Big((1+\bu^2)(L \psi_{\lambda,i})^2 +
    (1 + u^2)(\bL \psi_{\lambda,i})^2 \notag\\
    &\ \ \ \ \ \ \ \ \ \ \ \ \ \ + (1 + \bu^2 + u^2)(1-\frac{2M}{r})
    (\frac{\lambda^2}{r^2} + \frac{M}{r^3}
    )\cdot (\psi_{\lambda,i})^2 \Big)\ dr^* \ . \notag
\end{align}
Now, bringing the sum under the integral sign in the two terms on
the right hand side of \eqref{big_guy} above and then using the
$L^1$--$L^\infty$ H\"older inequality yields the bound:
\begin{multline}
        \sup_{0\leqslant s \leqslant t}\ \underline{E}(\phi(s)) \ \
        \
        \lesssim \ \ \ \underline{E}((1-\Delta_{sph})^\frac{1}{2}\phi(0)) \ + \
        \sup_{0\leqslant s\leqslant t}\ \big[ \underline{E}
        (\phi(s)) +
        E((1-\Delta_{sph})^\frac{1}{2}\phi(s)) \big]^\frac{1}{2}\\
        \cdot \int_0^t\ \lp{(1 + |\bu| + |u|)
    (1 - \frac{2M}{r} )\, r\cdot
    (\sqrt{1-\Delta_{sph}}G)(s)}{L^2(dr^*d\omega)}\ ds \ ,
    \label{little_guy}
\end{multline}
where the usual energy is given by:
\begin{equation}
        E(\phi(s)) \ = \ \int_{\RR\times\mathbb{S}^2 \times\{s\}}\
        \left(|\nabla_{t,r^*} (r\phi)|^2 + (1-\frac{2M}{r})
        \big( \frac{|\snabla_\omega (r\phi)|^2}{r^2} + \frac{M}{r^3}
        (r\phi)^2\big)\right)\ dr^* d\omega \ . \notag
\end{equation}
The estimate \eqref{morawetz} now follows from \eqref{little_guy}
and taking an angular (momentum) derivative of the basic energy
estimate (see the next section for a proof):
\begin{equation}
        \sup_{0\leqslant s \leqslant t}\ E(\phi(s)) \ \lesssim \
        E(\phi(0)) \ + \
        \left(
	\int_0^t \ \lp{(1-\frac{2M}{r})rG}{L^2(dr^* d\omega)}\ ds \right)^2
	\ . \notag
\end{equation}

\ret

\section{Proof of the Main Estimate}
We now turn to the proof of the estimate \eqref{main_thm_red}. This
will be accomplished in a series of three steps, each of which
represents a tightening of the usual energy estimate. These are:
\begin{enumerate}
        \item Usual conservation of energy.
    \item Weak local decay of energy.
    \item Strong uniform local decay of energy.
\end{enumerate}
Steps (1) and (3) involve a direct use of the energy-momentum tensor
identities recorded in the previous section applied to various
vector-fields $X$ which are associated with the various types of
decay as just listed. To prove item (2) above we use a 
Soffer--Morawetz type
multiplier similar to what was done in \cite{SB_1}--\cite{SB_2}.


\subsection*{Step 1: Conservation of energy}
This is well known. In the current setup it comes from setting
$X=T=\partial_t$  in \eqref{momentum_density}. Because $Q(x)$ is
time independent and since $T$ is Killing we end up with an
essentially divergence free quantity:
\begin{equation}
        \partial^\alpha {}^{(T)}P_\alpha \ = \
    H\cdot\partial_t(\psi) \ . \notag
\end{equation}
Integrating this over a time slab we arrive at the
energy estimate:
\begin{multline}
        \int_{\RR\times\{t\}}\ \left((L\psi)^2 + (\bL\psi)^2
    + Q(x)\cdot(\psi)^2\right)\ dx \\
    \lesssim \ \int_0^t\ \lp{H(s)}{L^2(dx)}\cdot
    \lp{\partial_t\psi(s)}{L^2(dx)}
    \ ds \ + \
    \int_{\RR\times\{0\}}\ \left((L\psi)^2 + (\bL\psi)^2
    + Q(x)\cdot(\psi)^2\right)\ dx
    \ , \label{basic_energy}
\end{multline}
where the implicit constant is fixed and does not depend on
$Q$ (it is easy to calculate but we
will not bother).\ret


\subsection*{Step 2: Weak Local Decay of Energy}
In this subsection, we prove that the local $L^1$ norm of the
quantity $Q\cdot(\psi)^2$ decays sufficiently fast in an average
sense. Our bound will be rather weak in that we allow the right hand
side of the estimate to grow like $\lambda t$. However, this weak
bound will be precisely what we need in the next subsection when we
prove the strong uniform local decay of energy. What we propose to
show is the following estimate for integers $1\leqslant N$:
\begin{multline}
        \int_0^t\int_{\mathcal{B}_1}\ (1+s) Q\cdot(\psi)^2 \ dx\, ds
    \ - \ \int_0^t\int_{\RR\setminus \mathcal{B}_1}\ (1+s)
    \chi_1(\frac{10x}{1+s})\ \sgn(x) \partial_x(Q) \cdot(\psi)^2\ dx\, ds
    \\ \\
    \lesssim \ \sup_{0\leqslant s \leqslant t}
    N^{-1}\underline{E}(s) \ + \ N (1 + \lambda^2) E(0)\\
    + N \int_0^t\int_{\RR}\
    \lp{(1 + s)(1+\lambda)H(s)}{L^2(dx)}
    \cdot\lp{(1 + \lambda)(|\nabla_{t,x}\psi| +
    Q^\frac{1}{2}\cdot|\psi|)\, (s)}{L^2(dx)}\ ds
    \ . \label{weak_local_decay}
\end{multline}
Here the implicit constant depends \emph{only} on the constants $C$
and the lengths of the interval $\mathcal{B}_1$ from lines
\eqref{rep1}--\eqref{rep6}. We are defining $\underline{E}$ as the
Morawetz type energy from line \eqref{mor_energy} above. Finally,
$\chi_1$ denotes a smooth bump adapted to the interval $[-1,1]$, and
$E$ denotes the basic energy from line \eqref{basic_energy} above.\\

In our proof of \eqref{weak_local_decay} it will be convenient for
us to make the assumption that the local bound \eqref{rep5}
drastically improves if we restrict to very small sets containing
$x=0$. In particular, we will assume that:
\begin{align}
        \frac{\epsilon^{-4}}{C}\cdot x^2 \ &\leqslant \
        -x\partial_x (Q)(x) \ , &|x| \ \leqslant \ c_{\mathcal{B}_1}
        \epsilon \ ,
        \label{crucial_ass}
\end{align}
for some sufficiently small parameter $\epsilon$ which will be
chosen in a moment. It is crucial for us to point out here that our
choice of $\epsilon$ will only be dictated by the constant $C$ and
the size of $\mathcal{B}_1$, and will not depend on any other
property of $Q$. Also, it is immediate that the assumption
\eqref{crucial_ass} in fact involves no loss of generality. This is
because the equation \eqref{1d_wave} rescales as follows:
\begin{align}
        \psi(t,x) \ &\rightsquigarrow \ \psi(\epsilon^{-1}t,\epsilon^{-1}x) \ ,
        &Q(x) \ &\rightsquigarrow \
        \epsilon^{-2}Q(\epsilon^{-1}x) \ ,
        &H(x) \ &\rightsquigarrow \
        \epsilon^{-2}H(\epsilon^{-1}t,\epsilon^{-1}x) \ . \notag
\end{align}
Notice that the conditions \eqref{rep1}--\eqref{rep6} adapt to the
rescaled situation in obvious ways. In particular one has
\eqref{crucial_ass} on the set $\widetilde{\mathcal{B}}_1
=\epsilon\cdot \mathcal{B}_1$. For the rest of this subsection we
will work in the rescaled situation where we assume all of
\eqref{rep1}--\eqref{rep6} as well as \eqref{crucial_ass}. Of course
once one has \eqref{weak_local_decay} in this rescaled situation,
one can recover the same bound for the original potential $Q$ by
scaling back. This will create constants which depend on $\epsilon$,
but  we choose this parameter  \emph{only} to overcome two things.
The first is the possibly large constant $C$ on the right hand side
of \eqref{rep5} (which is actually only a problem when $\lambda$ is
small). The second is the fact that the original $\mathcal{B}_1$ may
be small, so the constant $c_{\mathcal{B}_1}$ on line
\eqref{crucial_ass} where our improved bound holds is also small. Of
course both $C$ and $c_{\mathcal{B}_1}$ are fixed no matter how much
we rescale, so these can be made up for by taking $\epsilon$
sufficiently small. The main thing to keep in mind here is that our
rescaling  will never create constants in our estimates which depend
in other ways on the shape of $Q$, other than the original
assumptions we have made \eqref{rep1}--\eqref{rep6}.\\

To prove \eqref{weak_local_decay}, we use the following growth
multiplier of Soffer--Morawetz type:
\begin{equation}
        A(s,x)\psi \ = \
        (1+s) \chi_1(\frac{10x}{1+s}) \left[
        \varphi\partial_x \psi + \partial_x(\varphi \psi)
        \right] \ , \notag
\end{equation}
where $\varphi$ is defined as follows:
\begin{equation}
        \varphi(x) \ = \ \int_{0}^x\ \frac{1}{(1 + |y|)^k}\ dy \ ,
        \label{weight_function}
\end{equation}
where $1 <  k$ is a fixed constant. In practice the smaller the
value of $k$ the more favorable the estimates, so the reader may
assume that $k=2$. However, we will do all of our calculations in
the general case so the overall structure is more apparent. The
estimate \eqref{weak_local_decay} will  follow from the usual
procedure of directly calculating the integral:
\begin{align}
        I \ &= \ - \  \int_0^t \int_{\RR}\ H\cdot A(s,x)\psi\ dx\, ds \ ,
        \notag\\
            &= \ \int_0^t \int_{\RR}\
    \left(\partial_t^2\psi - \partial_x^2\psi\right)\cdot
    A(s,x)\psi\ dx\, ds \ + \ 
    \int_0^t \int_{\RR}\ Q\psi\cdot
    A(s,x)\psi\ dx\, ds \ , \notag\\
    &= \ I_1 + I_2 \ . \label{the_splitting}
\end{align}
and then using a Poincare type lemma near the critical point of
$Q(x)$ to get rid of the factor $-x\partial_x(Q)$ in favor of $Q$.
We now compute the terms $I_i$ separately and in order. The first
term $I_1$ is the pure error. We first integrate by parts with
respect to $\partial_t$ which yields the identity:
\begin{multline}
         I_1 \ = \  - \int_0^t \int_{\RR}\ \partial_t\psi\cdot
     A(s)(\partial_t\psi) \ dx\, ds
      - \int_0^t \int_{\RR}\ \partial_t\psi\cdot
     \dot{A}(s)(\psi) \ dx\, ds \\
     + \int_{\RR\times \{t\}}\ \partial_t\psi\cdot
     A(t)(\psi) \ dx -
     \int_{\RR\times \{0\}}\ \partial_t\psi\cdot
     A(0)(\psi) \ dx \ . \label{I_1_first}
\end{multline}
Here the operator $\dot{A}(s)$ is given by:
\begin{equation}
        \dot{A}(s)\psi \ = \ \left[\chi_1(\frac{10x}{1+s})
    - \frac{10x}{1+s}\chi'_1(\frac{10x}{1+s})\right]\cdot
    \left(2\varphi\partial_x\psi + \frac{1}{(1+|x|)^k}\psi
    \right) \ . \notag
\end{equation}
Also, one has the adjoint formula:
\begin{equation}
        A^*(s)\psi \ = \ -A\psi - 20 \chi'_1(\frac{10x}{1+s})
        \varphi\cdot\psi \ . \notag
\end{equation}
Therefore, a bound for the absolute value of the right hand side of
\eqref{I_1_first} above is:
\begin{multline}
        |I_1| \ \lesssim \ \int_0^t \int_{\RR}\
    \widetilde{\chi}_1(\frac{10x}{1+s})\cdot
    \left((\partial_t\psi)^2 + (\partial_x\psi)^2
    + \frac{1}{(1 + |x|)^{2k}} (\psi)^2\right) \ dx\, ds \\
    + \sup_{0\leqslant s\leqslant t}\int_{\RR\times\{s\}}\
    (1 + s) \widetilde{\chi}_1(\frac{10x}{1+s})
    \cdot \left((\partial_t\psi)^2 + (\partial_x\psi)^2
    + \frac{1}{(1 + |x|)^{2k}} (\psi)^2\right) \ dx \ ,
    \label{first_I_1_bound}
\end{multline}
where $\widetilde{\chi}_1$ is another $[-1,1]$ adapted smooth bump.
To deal with the terms involving the inverse $|x|$ weight we use the
Poincare type estimate:
\begin{equation}
        \int_{-x_0}^{x_0}\ \frac{1}{(1 + |x|)^2}\ (\psi)^2\ dx
        \ \lesssim \ (\psi)^2(0) + \int_{-x_0}^{x_0}\
        (\partial_x\psi)^2\ dx \ . \label{1_Poincare}
\end{equation}
This follows at once from evaluation of the integral:
\begin{multline}
        \frac{1}{(1 + |x_0|)}\left( (\psi)^2(x_0) + (\psi)^2(-x_0)\right)
        \ - \ 2(\psi)^2(0)\\
        = \ \int_{-x_0}^{x_0}\ \sgn(x)\partial_x\left[
        \frac{1}{(1 + |x|)}\ (\psi)^2\right]\ dx \ , \notag
\end{multline}
and using the Cauchy-Schwartz inequality. Using now
\eqref{1_Poincare} and the condition \eqref{rep6} it is easy to bound:
\begin{equation}
        \int_{\RR}\
        \widetilde{\chi}_1(\frac{10x}{1+s})\cdot
        \frac{1}{(1 + |x|)^{2k}} (\psi)^2 \ dx
        \ \lesssim \
        \int_{\RR}\
        \widetilde{\chi}_1(\frac{5x}{1+s})\cdot
        \left((\partial_x\psi)^2
        + Q\cdot (\psi)^2\right) \ dx \ , \label{cut_poincare}
\end{equation}
by using a partition of unity on $\mathcal{B}_1$ and
$\RR\setminus\mathcal{B}_1$.\\

Our next step is to use the bound:
\begin{equation}
        \int_{\RR}\ (1+s)\widetilde{\chi}_1(\frac{5x}{1+s})\cdot 
        \left((\partial_t\psi)^2 + (\partial_x\psi)^2
        + Q\cdot (\psi)^2\right)\ dx  \ \lesssim \
        (1 + s)^{-1}\cdot \underline{E}(s) \ , \notag
\end{equation}
where the right hand side is the Morawetz type energy from line
\eqref{mor_energy} above. This and the bound \eqref{first_I_1_bound}
allows us to estimate:
\begin{align}
        |I_1| \ &\lesssim \
        \int_{ N(1+\lambda)}^t\ (1 + s)^{-2}\cdot \underline{E}(s)\ ds \
        + \
        \sup_{  N(1+\lambda) \leqslant s \leqslant t}\ (1+s)^{-1}\cdot
        \underline{E}(s) \notag\\
        &\ \ \ \ \ \ \ + \ \int_0^{ N(1+\lambda)}\  E(s)\ ds \
        \ + \
        \sup_{ 0\leqslant s\leqslant N(1+\lambda) }\ (1+s)\cdot
        E(s) \ , \notag\\
        &\lesssim \ \sup_{0\leqslant s \leqslant t}
        (N(1+\lambda))^{-1}\underline{E}(s) \ + \ 
	\sup_{0\leqslant s \leqslant t}
        N (1 + \lambda) E(s) \ . \notag
\end{align}
Using now the energy inequality \eqref{basic_energy} to deal with
the second term on the right hand side of this last line we arrive
at the bound:
\begin{equation}
        |I_1|\ \lesssim \ (1 +
        \lambda)^{-1}\hbox{(R.H.S.)}\eqref{weak_local_decay} \ .
        \label{final_I1_bound}
\end{equation}
In a moment we will need to multiply all of our estimates through by
the factor $(1 + \lambda)$, so \eqref{final_I1_bound} is of the
correct form.\\

Before moving on to the second integral on line \eqref{the_splitting}
above, we mention briefly how to take care of the first integral on
the right hand side immediately above that line. Applying the
Cauchy-Schwartz inequality we have the bound:
\begin{multline}
        \int_0^t\int_{\RR}\ |H(s)|\cdot |A(s)\psi|\ dx \, ds \\
        \lesssim \ \int_0^t\ \lp{(1 + s)H(s)}{L^2(dx)}
        \cdot\lp{(|\partial_x \psi| + \chi_1 \, (1 + |x|)^{-k}|\psi|)(s)
        }{L^2(dx)} \ ds \  . \label{first_int}
\end{multline}
Using now a Poincare type estimate of the form \eqref{cut_poincare}
to deal with the last term on the right hand side of
\eqref{first_int} easily yields:
\begin{equation}
        \hbox{(L.H.S.)}\eqref{first_int} \ \lesssim \
        (1 + \lambda)^{-2} \ \hbox{(R.H.S.)}\eqref{weak_local_decay}
        \ , \notag
\end{equation}
which is sufficient for our purposes.\\

Finally, we deal with the integral $I_2$ on line
\eqref{the_splitting}. After several integration by parts (this is
an essentially well known calculation) we arrive at the identity:
\begin{equation}
        I_2 \ = \ \sum_{j=1}^4 \ K_j \ , \notag
\end{equation}
where the integrals $K_i$ are:
\begin{align}
        K_1 \ &= \ \int_0^t\int_{\RR}\
        10 \chi'_1(\frac{10x}{1+s}) \left[
        \varphi\partial_x \psi + \varphi' \psi
        \right]\cdot \partial_x\psi \ dx \ , \notag \\
        K_2 \ &= \ -\, \int_0^t\int_{\RR}\
        5 \chi'_1(\frac{10x}{1+s})\, \varphi'' \ (\psi)^2
        \ dx \ , \notag \\
        K_3 \ &= \ -\, \int_0^t\int_{\RR}\
        10 \chi'_1(\frac{10x}{1+s})\, \varphi\ Q\cdot (\psi)^2
        \ dx \ , \notag \\
        K_4 \ &= \  \int_0^t\int_{\RR}\
        (1+s) \chi_1(\frac{10x}{1+s})\left[
        2\varphi'(\partial_x\psi )^2 - \varphi\partial_x( Q)\cdot(\psi)^2
        - \frac{1}{2}\varphi''' (\psi)^2
        \right] \ dx \ . \notag
\end{align}
Bounding the first three terms above is essentially the same as what
we have just done for the term $I_1$ above. One simply uses
Cauchy-Schwartz, the Poincare estimate \eqref{1_Poincare}, and the
definitions of the two energies $E$ and $\underline{E}$ to prove
that:
\begin{equation}
        |K_1| + |K_2| + |K_3| \ \lesssim \
        \sup_{0\leqslant s \leqslant t}
        (N(1+\lambda))^{-1}\underline{E}(s) \ + \ \sup_{0\leqslant s \leqslant t}
        N (1 + \lambda) E(s) \ . \label{easy_K_bound}
\end{equation}
Therefore, the heart of the matter now is to obtain a positive lower
bound for the quantity $K_4$ in such a way that we can estimate the
left hand side of \eqref{weak_local_decay}.\\

Before continuing with the proof, let us make a further
simplification. Without loss of generality we may assume that the
cutoff function $\chi_1$ is the square of yet another smooth cutoff
function, say $\widetilde{\chi}_1$. This allows us to replace $\psi$
by $\td{\chi}_1\psi$ in $K_4$ above modulo a term involving
$[\partial_x,\td{\chi}_1]= O(\frac{1}{1+s})$ which is also cutoff
where $|x|\leqslant 10^{-1} t$. It is clear that this will again be
of the form
$(1+\lambda)^{-1}\hbox{(R.H.S.)}\eqref{weak_local_decay}$, so
we can just tack this error on to \eqref{easy_K_bound} above.\\

Thus,  what we will need to show is that there exists a sufficiently
small constant $c$ such that the following reverse bound holds for
compactly supported functions $\psi$:
\begin{multline}
        \int_{\RR}\
        \left[
        2\varphi'(\partial_x\psi )^2 - \varphi\partial_x( Q)\cdot(\psi)^2
        - \frac{1}{2}\varphi''' (\psi)^2
        \right] \ dx\\
         \ \geqslant \ c\ \int_{\RR}\
       \left[
        \varphi'(\partial_x\psi )^2 -  \frac{1}{2}
    \varphi\partial_x( Q)\cdot(\psi)^2
        \right] \ dx \ . \label{main_comm}
\end{multline}
Once this is established, the bound \eqref{weak_local_decay} will
follow from combining the bounds \eqref{final_I1_bound} and
\eqref{easy_K_bound} with \eqref{main_comm} and the following
estimate which also holds for smooth compactly supported functions
$\psi$:
\begin{multline}
         \int_{\mathcal{B}_1}\  Q\cdot(\psi)^2 \ dx
        \ - \ \int_{\RR\setminus \mathcal{B}_1}\
        \sgn(x) \partial_x(Q) \cdot(\psi)^2\ dx \\
         \lesssim \ (1+\lambda)\ \int_{\RR}\
       \left[
        \varphi'(\partial_x\psi )^2 - \frac{1}{2} \varphi\partial_x( Q)\cdot(\psi)^2
        \right] \ dx \ . \label{gets_lambda}
\end{multline}\ret

We first prove \eqref{main_comm}. The overall strategy for this is
very simple. The main thing we will establish is that the form of
the weight function \eqref{weight_function} reduces everything to
having a ``good'' bound for the function $\psi$ at $x=0$ in terms of
the hand side of \eqref{main_comm}. This latter task is relatively
easy to accomplish because the assumption \eqref{crucial_ass}
essentially means that $-x\partial_x(Q)\sim \epsilon^{-1}\delta_0$,
where $\delta_0$ is the unit mass at the origin. This means that the
potential term on the right hand side of \eqref{main_comm} will give
us a bound on $\psi(0)$ with an $O(\epsilon^\frac{1}{2})$ constant.
The details of this procedure are as follows. We first compute:
\begin{align}
        0 \ &= \ \int_{\RR} \ \partial_x\left[\varphi''\,
        (\psi)^2\right]\ dx \ , \notag\\
        &= \ \int_{\RR} \ \varphi'''\,
        (\psi)^2\ dx \ + \
        2\, \int_{\RR} \ \varphi''\,
        \psi \partial_x\psi\ dx \ . \label{start}
\end{align}
It will now be useful to have the identities:
\begin{align}
        \varphi''(x) \ &= \ \frac{-k \cdot\sgn(x)}{(1+|x|)^{k+1}} \ ,
        &\varphi'''(x) \ &= \ -2k\delta_0 \ + \ \frac{k(k+1) }{(1+|x|)^{k+2}} \
        . \notag
\end{align}
Therefore, the right hand side of \eqref{start} and a
Cauchy-Schwartz gives us the bound:
\begin{align}
        &\int_{\RR} \ \frac{k(k+1)}{(1+|x|)^{k+2}}\,
        (\psi)^2\ dx \ , \notag\\
        \leqslant \
        &2\left(\frac{k}{k+1}\right)^\frac{1}{2}\left(\int_{\RR}\
        \frac{1}{(1+|x|)^{k}}\, (\partial_x\psi)^2\ dx
        \right)^\frac{1}{2}
        \cdot\left(\int_{\RR}\ \frac{k (k+1)}{(1+|x|)^{k+2}}
        (\psi)^2\ dx \right)^\frac{1}{2}
        + 2k(\psi)^2(0) \ , \notag\\
        = \ &2\left(\frac{k}{k+1}\right)^\frac{1}{2}
        A^\frac{1}{2}\cdot B^\frac{1}{2} + C \ . \notag
\end{align}
We may now assume without loss of generality that in this last bound
we have $C\leqslant B$, otherwise there is nothing to prove on line
\eqref{main_comm}. Therefore, dividing through by $B^\frac{1}{2}$
and squaring this last line we arrive at the bound:
\begin{equation}
        \left|\int_{\RR} \ \varphi'''\,
        (\psi)^2\ dx  \right| \ = \
        B-C \ \leqslant \ 4\frac{k}{k+1}A
        + 4 \left(\frac{k}{k+1}\right)^\frac{1}{2}
    A^\frac{1}{2}\cdot C^\frac{1}{2} \ . \label{main_thing}
\end{equation}
The dangerous term is now the second one on the right hand side
above. This needs to be controlled in terms of a sufficiently small
constant. In fact, we will show that it is $O(\epsilon^\frac{1}{2})$
times the (R.H.S.)\eqref{main_comm}, which implies that it may be
safely absorbed into half of the remaining portion of $A$ and a
small amount of the potential term on (R.H.S.)\eqref{main_comm}. The
bound which allows us to do this is the following:
\begin{equation}
        (\psi)^2(0) \ \lesssim \
        \epsilon\ \left(\int_{\RR}\
        \frac{1}{(1+|x|)^{k}}\, (\partial_x\psi)^2\ dx +
        \int_{\RR}\
        -\varphi \partial_x(Q)\cdot (\psi)^2\ dx\right) \ . \notag
\end{equation}
From the assumption \eqref{crucial_ass}, this last estimate follows
from:
\begin{equation}
        (\psi)^2(0) \ \lesssim \
        \epsilon\ \left(\int_{\RR}\
        \frac{1}{(1+|x|)^{k}}\, (\partial_x\psi)^2\ dx +
        \epsilon^{-4}\ \int_{\RR}\
        x^2\chi(\epsilon^{-1}x) \cdot (\psi)^2\ dx\right) \ ,
        \label{special_sob}
\end{equation}
where $\chi$ is some smooth $O(1)$ bump function whose support
depends on the size of the set $\mathcal{B}_1$ from line
\eqref{rep5}.  The estimate \eqref{special_sob} is essentially scale
invariant, so it suffices to show that:
\begin{equation}
        (\psi)^2(0) \ \lesssim \
        \int_{\RR}\
        \frac{1}{(1+\epsilon|x|)^{k}}\, (\partial_x\psi)^2\ dx +
        \ \int_{\RR}\
        x^2\chi(x) \cdot (\psi)^2\ dx \ .
        \label{rescaled_guy}
\end{equation}
This, in turn, follows from cutting things off and using the usual
Sobolev embedding once we have the bound:
\begin{equation}
        \lp{\td{\chi}^\frac{1}{2}\psi}{L^2}^2 \ \lesssim \
        \int_{\RR}\
        \frac{1}{(1+\epsilon|x|)^{k}}\, (\partial_x\psi)^2\ dx +
        \ \int_{\RR}\
        x^2\chi(x) \cdot (\psi)^2\ dx \ ,
        \label{final_sob}
\end{equation}
for some slightly smaller cutoff function $\td{\chi}$. This last
bound can be proved in two steps. We first show the estimate:
\begin{equation}
        \lp{|x|^\frac{1}{2}\, \td{\td{\chi}}^\frac{1}{2} \psi}{L^2}^2 \ \lesssim \
        \int_{\RR}\
        \frac{1}{(1+\epsilon|x|)^{k}}\, (\partial_x\psi)^2\ dx +
        \ \int_{\RR}\
        x^2\chi(x) \cdot (\psi)^2\ dx \ ,
        \label{inter_sob}
\end{equation}
for some intermediate cutoff $\td{\td{\chi}}$. This bound follows at
once from evaluating the integral:
\begin{equation}
        0 \ = \ \int_{\RR}\ \sgn(x)\, \partial_x\left[
        x^2\ \td{\td{\chi}} \cdot(\psi)^2
        \right]\ dx \ , \notag
\end{equation}
and using the Cauchy-Schwartz inequality to bound the error terms by
(R.H.S.)\eqref{inter_sob}. Having now established \eqref{inter_sob}
we can prove \eqref{final_sob} by applying the same procedure to the
integral:
\begin{equation}
        0 \ = \ \int_{\RR}\  \partial_x \left[
        x\ \td{\chi} \cdot(\psi)^2
        \right]\ dx \ . \notag
\end{equation}
This completes our proof of \eqref{special_sob}, and hence our
demonstration of the main commutator estimate \eqref{main_comm}.\\

Having now dealt with the bound \eqref{main_comm}, the only thing
left for us to do in our proof of \eqref{weak_local_decay} is to
show the bound \eqref{gets_lambda}. Notice that the bound for the
second term on the left hand side of that estimate follows at once
from the  fact that $\sgn(x)\lesssim \varphi(x)$ whenever
$x\in\RR\setminus \mathcal{B}_1$. Therefore, it remains to bound the
first term on the left hand side of \eqref{gets_lambda}. This is
where we pick up the extra factor of $(1+\lambda)$. The proof is
essentially identical to what was done to establish
\eqref{final_sob} above. Using the two conditions
\eqref{rep5}--\eqref{rep6}, it suffices to multiply through the
following estimate by the quantity $(1+\lambda)^2$:
\begin{equation}
         \lp{ \chi^\frac{1}{2}_{\mathcal{B}_1}\psi}{L^2}^2
        \ \lesssim \ \lp{ |x|
        {\chi}^\frac{1}{2}_{\mathcal{B}_1}
        \psi}{L^2}\cdot\lp{
        {\chi}^\frac{1}{2}_{\mathcal{B}_1}\partial_x
        \psi}{L^2}
        + \lp{ |x|
        \td{\chi}^\frac{1}{2}_{\mathcal{B}_1}
        \psi}{L^2}^2
        \ . \label{L2_est}
\end{equation}
Here the functions $\chi_{\mathcal{B}_1}$ and
$\td{\chi}_{\mathcal{B}_1}$ are cutoffs which are $\equiv 1$ on the
set $\mathcal{B}_1$ and which vanish outside of $2\mathcal{B}_1$.
The bound \eqref{L2_est} follows from evaluation of the integral:
\begin{equation}
         0 \ = \   \int_{\RR}\  \partial_x \left[
        x\ \chi_{\mathcal{B}_1} \cdot(\psi)^2
        \right]\ dx \ , \notag
\end{equation}
and using Cauchy-Schwartz as well as the bound
$|x\chi'_{\mathcal{B}_1}|\lesssim x^2 \td{\chi_{\mathcal{B}_1}}$ for
a suitable cutoff $\td{\chi}_{\mathcal{B}_1}$. We have now finished
our proof of the weak local energy decay estimate
\eqref{weak_local_decay}.\\

\begin{rem}
We note here that it is possible to prove \eqref{weak_local_decay}
without rescaling the potential $Q$ into the condition
\eqref{crucial_ass}. This can be accomplished by using the weight
function:
\begin{equation}
        \varphi_\epsilon(x) \ = \ \int_{0}^x\ \frac{1}
	{(1 + \epsilon|y|)^k}\ dy \ ,
        \notag
\end{equation}
in place  of \eqref{weight_function} above. This yields a small
factor in front of $|\psi|(0)$ when it appears in the
$C^\frac{1}{2}$ term on the right hand side of line
\eqref{main_thing} above, so one can proceed directly to the
estimate \eqref{rescaled_guy} to control things. We leave the
details to the interested reader.
\end{rem}\ret


\subsection*{Step 3: Strong Uniform Decay of Local Energy}

We are now ready to prove the main Morawetz estimate \eqref{1D_mor}.
With the assumptions \eqref{rep1}--\eqref{rep6} in hand, as well as
the weak local energy decay estimate \eqref{weak_local_decay}, this
becomes an essentially standard calculation. We will contract the
energy-momentum tensor \eqref{1d_em} with the conformal Killing
vector-field:
\begin{equation}
        K_0 \ = \ (t^2 + x^2)\partial_t + 2tx\partial_x \ = \
        \frac{1}{2}\bu^2 L + \frac{1}{2} u^2\bL \ . \notag
\end{equation}
The deformation tensor of this is computed to be:
\begin{equation}
        {}^{(K_0)}\pi \ = \ 4t g \ . \notag
\end{equation}
Therefore, we may form the momentum density ${}^{(K_0)}P_\alpha =
T_{\alpha\beta} K_0^\beta$ and from line \eqref{div_iden} we compute
the divergence:
\begin{equation}
        \partial^\alpha \, {}^{(K_0)}P_\alpha \ = \
        -tx\partial_x(Q)\cdot(\psi)^2 - 2tQ\cdot(\psi)^2 +
        K_0(\psi)\cdot H \ . \label{K0_div_iden}
\end{equation}
By simply integrating this last line over various time slabs of the
form $0\leqslant s \leqslant t$ and using the Cauchy-Schwartz
inequality we arrive at the bound:
\begin{multline}
        \sup_{0\leqslant s \leqslant t}\
        {}^{(K_0)}P_0 (s)\ \leqslant \
        \int_0^t\int_{\RR}\ |H(s)|\cdot|K_0(\psi)(s)|\ dx ds
        + {}^{(K_0)}P_0(0)\\
        + \int_0^t\int_{\RR}\ \left[
        sx\partial_x(Q)\cdot(\psi)^2 + 2sQ\cdot(\psi)^2
        \right]\ dx\, ds \ . \label{mor_div_iden}
\end{multline}
Using now the identity:
\begin{equation}
        {}^{(K_0)}P_0 \ = \ \frac{1}{4}\bu^2(L\psi)^2 +
        \frac{1}{4}u^2(\bL\psi)^2 + \frac{1}{4}(\bu^2 +
        u^2)Q\cdot(\psi)^2 \ , \notag
\end{equation}
we see that \eqref{mor_div_iden} in conjunction with the energy
estimate \eqref{basic_energy} implies the bound:
\begin{multline}
        \sup_{0\leqslant s \leqslant t}\
        \underline{E} (s)\ \lesssim \
        \int_0^t\int_{\RR}\ \lp{(1 + |\bu| + |u|)H(s)}{L^2(dx)}\cdot
        \underline{E}^\frac{1}{2}(s)\ ds
        + \underline{E}(0)\\
        + \int_0^t\int_{\RR}\ \left[
        sx\partial_x(Q)\cdot(\psi)^2 + 2sQ\cdot(\psi)^2
        \right]\ dx\, ds \ . \label{almost_there}
\end{multline}
The last thing we need to do here is to bound the last term on the
right hand side of the previous expression. We will show the bound:
\begin{multline}
        \int_0^t\int_{\RR}\ \left[
        sx\partial_x(Q)\cdot(\psi)^2 + 2sQ\cdot(\psi)^2
        \right]\ dx\, ds \ \lesssim \
        \sup_{0\leqslant s \leqslant t}
        N^{-1}\underline{E}(s) \ + \ N (1 + \lambda^2) E(0)\\
        + N \int_0^t\int_{\RR}\
        \lp{(1 + s)(1+\lambda)H(s)}{L^2(dx)}
        \cdot\lp{(1 + \lambda)(|\nabla_{t,x}\psi| +
        Q^\frac{1}{2}\cdot|\psi|)\, (s)}{L^2(dx)}\ ds \ . \label{last_step}
\end{multline}
where the implicit constant is independent of the large parameter
$N$. Notice that this bound substituted into \eqref{almost_there}
immediately implies \eqref{1D_mor} for sufficiently large $N$.\\

To prove \eqref{last_step} we will chop the left hand side up into
three pieces. The first is the ``bad'' set $\mathcal{B}_1$. This is
where most of the positivity of (L.H.S.)\eqref{last_step} can be
found. The second set is where $x\notin\mathcal{B}_1$ and
$|x|\leqslant \frac{1}{20}t$. Here we use the strongly repulsive
condition \eqref{rep3}. Finally, in the exterior of the influence of
the potential when $t\lesssim |x|$ we can simply integrate things
using the homogeneity bound \eqref{rep4}. The details of this
procedure are as follows:\\

On the set $\mathcal{B}_1$ we use the repulsive condition
\eqref{rep2} and the first term  on the left hand side of
\eqref{weak_local_decay} to bound:
\begin{align}
        &\int_0^t\int_{\mathcal{B}_1}\ \left[
        sx\partial_x(Q)\cdot(\psi)^2 + 2sQ\cdot(\psi)^2
        \right]\ dx\, ds \  , \notag\\
        \lesssim \
        &\int_0^t\int_{\mathcal{B}_1}\ (1+s) Q\cdot(\psi)^2
        \ dx\, ds \ , \notag\\
        \lesssim \ &\hbox{(R.H.S.)}\eqref{last_step} \ . \notag
\end{align}
Next, we work in the set $\RR\setminus \mathcal{B}_1$ but cutoff
according to how large $t$ is. Here we make use of the condition 
\eqref{rep3}:
\begin{align}
        &\int_0^t\int_{\RR\setminus \mathcal{B}_1}\ \chi_1(
        \frac{10x}{1+s})\left[
        sx\partial_x(Q)\cdot(\psi)^2 + 2sQ\cdot(\psi)^2
        \right]\ dx\, ds \ , \notag\\
        \lesssim \ &-\ \int_0^t\int_{\RR\setminus \mathcal{B}_1}\ (1+s)
        \chi_1( \frac{10x}{1+s})\sgn(x)\partial_x(Q)\cdot(\psi)^2
        \ dx\, ds \ , \notag\\
        \lesssim \ &\hbox{(R.H.S.)}\eqref{last_step} \ . \notag
\end{align}
Finally, in the exterior where $t \lesssim |x|$ we use the condition
\eqref{rep3} and the following bound which holds for parameters $N$
such that $|\mathcal{B}_2|\leqslant N$ (where the implicit constant
of course depends on $|\mathcal{B}_2|$):
\begin{align}
        &\int_0^t\int_{\RR}\ \left( 1 - \chi_1(
        \frac{10x}{1+s})\right) \left[
        sx\partial_x(Q)\cdot(\psi)^2 + 2sQ\cdot(\psi)^2
        \right]\ dx\, ds \ , \notag\\
        \lesssim \ &\int_0^{20|\mathcal{B}_2|}\int_{\RR}\  2sQ\cdot(\psi)^2
        \ dx\, ds \ + \
        \int_{0}^{N}\int_{\RR}\  Q\cdot(\psi)^2
        \ dx\, ds \ + \
        \int_{N}^t\int_{\RR}\  Q\cdot(\psi)^2
        \ dx\, ds \ \notag\\
        \lesssim \ &N \sup_{0\leqslant s \leqslant t}E(s) \ + \
        N^{-1} \sup_{0\leqslant s \leqslant t}\underline{E}(s)
        \ , \notag\\
        \lesssim \ &\hbox{(R.H.S.)}\eqref{last_step} \ . \notag
\end{align}
This completes our demonstration of \eqref{last_step}, and hence our
proof the main estimate \eqref{1D_mor}.

\ret

\section{Scattering for Small Amplitude
Non-Linear Scalar Fields}We will be brief here and leave many of the
details to the reader. The main result of this section is the
following:

\begin{thm}[Scattering for Scalar Fields]\label{NLW_thm}
Consider the Cauchy problem:
\begin{align}
        \Box_g\phi \ &= \ \lambda |\phi|^p\phi \ ,
        &\phi(0) \ &= \ f \ ,
        &\partial_t\phi(0) \ &= \ g \ , \label{cp}
\end{align}
for compactly supported functions $(f,g)$. Define the regularity
space:
\begin{equation}
        \lp{\phi}{\mathcal{H}_\Omega^k}^2 \ = \
        \underline{E}
        \left((1-\Delta_{sph})^\frac{k}{2}\phi\right) \ , \notag
\end{equation}
where $\underline{E}$ is the Morawetz type energy from line
\eqref{morawetz}. Then if $2<p$, there exists a universal set of
positive constants $\mathcal{E}$ and $C$ depending only on $p$ such
that if:
\begin{equation}
        \lp{\phi(0)}{\mathcal{H}^3_\Omega} \ \leqslant \ \mathcal{E} \ , \notag
\end{equation}
then a unique solution to the problem \eqref{cp} exists for all
values of the variable $t$ and it obeys the bound:
\begin{equation}
        \lp{\phi(t)}{\mathcal{H}^2_\Omega} \ \leqslant \ C\mathcal{E} \
        . \notag
\end{equation}
In particular, one has the following uniform point-wise bounds:
\begin{equation}
        |\phi| \ \lesssim \ \mathcal{E}\ r^{-1}\cdot \min\{ \ 1  \ , \
        \big|t-|r^*|\big|^{-\frac{1}{2}}\  \} \ . \label{Linfty_est}
\end{equation}
\end{thm}\ret

\noindent A previous  result of this type was recently obtained by
Dafermos and Rodnianski in the case of spherical symmetry and powers
$3< p$ (see \cite{DR_2}).\\

To prove Theorem \eqref{NLW_thm} we need four ingredients. The first
is the Morawetz estimate \eqref{morawetz}. The second is a Poincare
type estimate which will allow us to control the $L^2$ norm of our
function in terms of the energy $\underline{E}$. The third is a
paraproduct bound which allows us to concentrate all of our angular
derivatives on a single term of the non-linearity $\lambda
|\phi|^p\phi$. And the final is a global Sobolev inequality which
will give us the bound \eqref{Linfty_est} in terms of our energy
space $\mathcal{H}_\Omega^2$. We now state the last three of these
in order:\\

\begin{lem}[Poincare type estimate for the weights $\bu$ and $u$]
Let $\psi$ be a function of the variable $r^*$, and define the
weights $\bu$ and $u$ as on lines \eqref{defs}. Then the following
estimate holds:
\begin{equation}
        \int_{\RR}\ (\psi)^2\ dr^* \ \lesssim \
        \int_{\RR}\ \left( \bu^2 (L\psi)^2 +
        u^2(\bL\psi)^2 + \frac{(1+
	  \bu^2 + u^2)(1-\frac{2M}{r})}{r^3}(\psi)^2 \right)
        \ dr^* \ . \label{the_poincare}
\end{equation}
\end{lem}\ret

\begin{lem}[Paraproduct bounds]
On the sphere $\mathbb{S}^2$ the following estimates holds:
\begin{equation}
        \lp{(1-\Delta_{sph})^\frac{k}{2} (|F|^p F) }{L^2(\mathbb{S}^2)}
        \ \lesssim \ \lp{F}{L^\infty(\mathbb{S}^2)}^p
        \cdot \lp{(1-\Delta_{sph})^\frac{k}{2} F}{L^2(\mathbb{S}^2)}
        \ , \label{basic_para}
\end{equation}
for all test functions $F$ and integers $0\leqslant k \leqslant
p+1$.
\end{lem}\ret

\begin{lem}[A global Sobolev inequality]
Let $\phi$ be a function of the variables $(r^*,\omega)$. Then one
has the following global bounds:
\begin{equation}
        |r\phi| \ \lesssim \   \min\{ \ 1 \ ,
        \ r^\frac{1}{2}(1-\frac{2M}{r})^{-\frac{1}{4}}
        \big|t-|r^*|\big|^{-1} \ , \
        \big|t-|r^*|\big|^{-\frac{1}{2}}\  \}\cdot
        \lp{\phi}{\mathcal{H}_\Omega^2} \ . \label{global_sob}
\end{equation}
\end{lem}\ret

\noindent We now give short proofs of these three Lemmas:\\

\begin{proof}[Proof of estimate \eqref{the_poincare}]
The proof will follow from cutting the function $\psi$ into three
pieces. We write:
\begin{equation}
        \psi \ = \ \chi_{r^* < -1}\psi + \chi_{-1 < r^* < 1}\psi
        + \chi_{r^* > 1}\psi \ , \notag
\end{equation}
where the $\chi$ form a smooth partition of unity such that the
$\partial_{r^*} \chi_{r^* < \pm 1}$ are supported on the interval
$[-2,2]$. For the left hand portion we compute that:
\begin{align}
        0 \ &= \ \int_{\RR}\
        \partial_{r^*}\left[ (t + r^*)\chi_{r^* < -1}(\psi)^2
        \right]\ dr^* \ , \notag\\
        &= \ \int_{\RR}\
        \chi_{r^* < -1}(\psi)^2
        \ dr^* \ + \ \int_{\RR}\
        (t+r^*)\chi'_{r^* < -1}(\psi)^2 \ dr^* \ + 2 \int_{\RR}\ (t+r^*)
       \chi_{r^* < -1}\psi \partial_{r^*}\psi \ dr^* \ . \notag
\end{align}
Collecting terms and applying the Cauchy-Schwartz inequality we
arrive at the bound:
\begin{multline}
        \int_{\RR}\
        \chi_{r^* < -1}(\psi)^2
        \ dr^*  \ \lesssim \
        \int_{\RR}\  \frac{(1+
	\bu^2 + u^2)(1-\frac{2M}{r})}{r^3}(\psi)^2 \ dr^*\\
        + \left(
        \int_{\RR}\ \chi_{r^* < -1}(\psi)^2 \ dr^*
        \right)^\frac{1}{2}\cdot\left(
        \int_{\RR}\ \bu^2 \ \chi_{r^* < -1}(\partial_{r^*}\psi)^2\  dr^*
        \right)^\frac{1}{2} \ . \notag
\end{multline}
This easily proves \eqref{the_poincare} for the left hand portion of
things because here one has the bound:
\begin{equation}
        \bu^2 \ \chi_{r^* < -1}(\partial_{r^*}\psi)^2 \ \lesssim \
        \bu^2(L\psi)^2 + u^2(\bL\psi)^2 \ . \label{bu_bound}
\end{equation}
The proof of \eqref{the_poincare} for the right hand term $\chi_{r^*
> 1}\psi$ follows from an identical argument. We leave this to the
reader.
\end{proof}\ret

\begin{proof}[Proof of estimate \eqref{basic_para}]
Because $k$ is chosen to be an integer, it is easy to reduce this
estimate to the case of Euclidean space. First of all, note that one has
the following bound:
\begin{equation}
    \lp{(1-\Delta_{sph})^\frac{k}{2} (|F|^p F) }{L^2(\mathbb{S}^2)}
    \ \lesssim \ \sum_{|\alpha|\leqslant k}
    \lp{\Omega_{ij}^\alpha (|F|^p F)}{L^2(\mathbb{S}^2)} \ ,
    \notag
\end{equation}
where $\alpha$ is a multiindexing of the rotation generators
$\{\Omega_{ij}\}$. Via a partition of unity, the desired bound
easily reduces to proving that for $k\in\mathbb{N}$ one has:
\begin{align}
        \lp{|F|^p F}{H^k} \ &\lesssim \ \lp{F}{L^\infty}^p
        \cdot\lp{F}{H^k} \  &0\leqslant k\leqslant p+1
        \ . \notag
\end{align}
Estimates of this type are well known and easy to prove (see e.g.
\cite{T_tools}).
\end{proof}\ret

\begin{proof}[Proof of estimate \eqref{global_sob}]
We first make a preliminary reduction. Because we are including two
angular (momentum) derivatives in the norm $\mathcal{H}_\Omega^2$,
via the Sobolev embedding on the sphere $\mathbb{S}^2$ it suffices
to prove the following global Sobolev estimate for functions $\psi$
of the variable $r^*$:
\begin{equation}
         |\psi| \ \lesssim \   \min\{ \ 1 \ ,
        \ r^\frac{1}{2}(1-\frac{2M}{r})^{-\frac{1}{4}}
        \big|t-|r^*|\big|^{-1} \ , \
        \big|t-|r^*|\big|^{-\frac{1}{2}}\  \}\cdot
        \lp{\phi}{\mathcal{H}} \ , \label{psi_linfty}
\end{equation}
where $\mathcal{H}$ is the Hilbert Space:
\begin{equation}
        \lp{\psi}{\mathcal{H}}^2\ = \ \int_{\RR}\
        \left[   (1+ \bu^2) (L\psi)^2 +
        (1+u^2)(\bL\psi)^2 + (1 + \bu^2 + u^2)\frac{(1-\frac{2M}{r})}{r^3}(\psi)^2
        \right]\ dr^* \ . \notag
\end{equation}
Using now the Poincare estimate \eqref{the_poincare} and bounds of
the form \eqref{bu_bound} we see that we have:
\begin{equation}
        \int_{\RR}\
        \left[   (1 + (t-|r^*|)^2)(\partial_r^*\psi)^2
        \ +\  (\psi)^2 \ +\
         (1 + \bu^2 + u^2)\frac{(1-\frac{2M}{r})}{r^3}(\psi)^2
        \right]\ dr^* \ \lesssim \
        \lp{\psi}{\mathcal{H}}^2\ . \label{main_energy}
\end{equation}
The first two terms on the left hand side of this last expression
are enough to get the basic weight in the estimate
\eqref{psi_linfty}. Specifically, one has the bound:
\begin{equation}
        (1 + (t-|r^*|)^2)^\frac{1}{2} (\psi)^2 \ \lesssim \
        \int_{\RR}\
        \left[   (1 + (t-|r^*|)^2)(\partial_r^*\psi)^2
        \ +\  (\psi)^2 \right]\ dr^* \ . \notag
\end{equation}
The proof of this kind of estimate is completely standard and
reduces to the usual Sobolev bound after decomposing things via
cutoffs on intervals of dyadic sizes according to the value of the
weight $(1 + (t-|r^*|)^2)^\frac{1}{2}$. See \cite{LS_MKG} for
details on this procedure. \\

It remains for us to prove the estimate \eqref{psi_linfty} for the
more refined weight. It is clear from the form of the energy
\eqref{main_energy} that this bound holds when $r^*\in[-1,1]$.
Therefore we only need to establish it for the cases $r^*<-1$ and
$r* > 1$. We will do this separately by directly evaluating some
weighted integrals similar to what we have done many times now. For
values $r^*_0< -1$ we compute that (note that the weight $r$ is
essentially a constant here so we can safely disregard it):
\begin{multline}
        (1-\frac{2M}{r})^\frac{1}{2}(t+r^*_0)^2(\psi)^2(r^*_0) \\
        = \ -\ \int_{r^*_0}^0\
        \partial_{r^*}\left[ (1-\frac{2M}{r})^\frac{1}{2}(t+r^*)^2 (\psi)^2
        \right]\ dr^* \ + \
        (1-\frac{2M}{r})^\frac{1}{2}\ t^2 \ (\psi)^2(0) \ . \notag
\end{multline}
By collecting terms and using the Cauchy--Schwartz inequality this
gives us the bound:
\begin{multline}
        (1-\frac{2M}{r})^\frac{1}{2}(t+r^*_0)^2(\psi)^2(r^*_0) \
        \lesssim \
        (1-\frac{2M}{r})^\frac{1}{2}\ t^2 \ (\psi)^2(0)\\
        + \ \lp{(1-\frac{2M}{r})^\frac{1}{2}(t+r^*)\psi}{L^2([0,r^*_0])}
        \cdot\lp{(t+r^*) \partial_{r^*} \psi}{L^2([0,r^*_0])} \  \ . \notag
\end{multline}
the crucial thing to notice here is that when the derivative falls
on the weight $(1-\frac{2M}{r})(t+r^*)^2$ it creates a
\emph{positive} term which can be collected with the left hand side.
This last line gives us the desired bound because all the $L^2$ type
norms are covered by the energy $\mathcal{H}$, and so is the bound
for $\psi$ at the
origin.\\

To wrap things up here we need to prove the refined bound in
\eqref{psi_linfty} for the region where $1\leqslant r^*$. Notice
that it suffices to do this for $r^*\leqslant t$ because otherwise
the simpler weight on the right hand side of \eqref{psi_linfty} is
favorable. Also, in this case we may safely ignore the weight
$(1-\frac{2M}{r})$ because it is essentially a non-zero constant.
The desired bound will drop out from computing the following
integral for fixed points $0\leqslant r^*_0\leqslant t$:
\begin{equation}
        r^{-1}(t-r^*_0)^2(\psi)^2(r^*_0) \ = \
        \int_0^{r^*_0}\ \partial_{r^*}\left[ r^{-1}(t-r^*)^2(\psi)^2
        \right]\ dr^* \ + \ r^{-1}\ t^2 \ (\psi)^2(0) \ .
        \notag
\end{equation}
Computing the integral on the right hand side of this last
expression, throwing away \emph{only} the term which results when
the derivative falls on the weight $(t-r^*)^2$, collecting the terms
of like sign to the left hand side, and applying the Cauchy-Schwartz
inequality we arrive at the bound:
\begin{multline}
        r^{-1}(t-r^*_0)^2(\psi)^2(r^*_0) \
        + \ \lp{ r^{-1}(t-r^*)\psi}{L^2([0,r^*_0])}^2
        \lesssim \ r^{-1}\ t^2 \
        (\psi)^2(0) \\
        + \ \lp{ r^{-1}(t-r^*)\psi}{L^2([0,r^*_0])}
        \cdot\lp{ (t-r^*) \partial_{r^*}\psi}{L^2([0,r^*_0])} \ .
        \notag
\end{multline}
From this one easily derives the bound:
\begin{equation}
        r^{-1}(t-r^*_0)^2(\psi)^2(r^*_0) \ \lesssim \
        r^{-1}\ t^2 \ (\psi)^2(0) \ + \ \lp{ (t-r^*) 
	\partial_{r^*}\psi}{L^2([0,r^*_0])}
        \ , \notag
\end{equation}
which is itself bounded by the energy $\mathcal{H}$. This completes
our proof of the estimate \eqref{psi_linfty}, and hence our
demonstration of the global Sobolev bound \eqref{global_sob}.
\end{proof}\ret\ret

\begin{proof}[Proof of Theorem \ref{NLW_thm}]
We now use the previous three lemmas and the main decay estimate
\eqref{morawetz} to prove the global regularity result. This will
follow by bootstrapping the usual local existence theorem. We will not
state or prove this local result here because even in this 
context it is an elementary
application of Picard iteration and energy estimates (in fact, one can
use the estimates developed here to set up a global Picard
iteration because the precise
structure of the non-linearity does not need to be preserved).
Now, our theorem will follow if we can show that for each fixed time
$t$ up to which we have existence the weak bound
$\sup_{0\leqslant s\leqslant t}\lp{\phi(s)}{\mathcal{H}_\Omega^2}
\leqslant 2C\mathcal{E}$ implies the
stronger bound $\sup_{0\leqslant s\leqslant t}\lp{\phi(s)}
{\mathcal{H}_\Omega^2}\leqslant C\mathcal{E}$. From
the estimate \eqref{morawetz}, we see that we can provide this as long
as we can show the bound:
\begin{multline}
        \int_0^t \ \lp{(1 + |\bu| + |u|)
    (1 - \frac{2M}{r} )\, r\cdot
    (\sqrt{1-\Delta_{sph}})^\frac{3}{2}(|\phi|^p\phi)\ (s)}
    {L^2(dr^*d\omega)}\ ds \\ \lesssim \
    \sup_{0\leqslant s\leqslant t}
    \lp{\phi(s)}{\mathcal{H}_\Omega^2}^{p+1} \ . \label{need_to_bound}
\end{multline}
To compute the integral on the left hand side, we use the paraproduct
estimate \eqref{basic_para} which gives us the bound:
\begin{equation}
        \lp{(\sqrt{1-\Delta_{sph}})^\frac{3}{2}(|\phi|^p\phi)\ (s)}
    {L^2(d\omega)} \ \lesssim  \
    \lp{\phi(s)}{L^\infty_\omega}^p
    \cdot\lp{(\sqrt{1-\Delta_{sph}})^\frac{3}{2}\phi\
    (s)}{L^2(d\omega)} \ . \notag
\end{equation}
Using now the definition of the energy $\mathcal{H}^2_\Omega$ and the
Poincare estimate \eqref{the_poincare} to control the $L^2$ norm of
the zero harmonic of $\phi$, we see that we can bound:
\begin{equation}
        \hbox{(L.H.S.)}
    \eqref{need_to_bound}
    \ \lesssim \
    \sup_{0\leqslant s\leqslant t}
    \lp{\phi(s)}{\mathcal{H}_\Omega^2}
    \cdot
    \int_0^t\
    \lp{(1 + |\bu| + |u|)
    (1 - \frac{2M}{r} )^\frac{1}{2} |\phi|^p \ (s)}
    {L^\infty_{(r^*,\omega)}}\ ds \ . \notag
\end{equation}
The claim \eqref{need_to_bound} will now follow once we can show the
fixed time bound:
\begin{equation}
        \lp{(1 + |\bu| + |u|)
    (1 - \frac{2M}{r} )^\frac{1}{2} |\phi|^p \ (s)}
    {L^\infty_{(r^*,\omega)}}
    \ \lesssim \
    (1+s)^{-1 -\frac{1}{2}(p-2)}\lp{\phi(s)}{\mathcal{H}_\Omega^2}^p \
    . \label{end_linfty_bound}
\end{equation}
This last bound follows easily from the global Sobolev estimate
\eqref{global_sob}. To see this, it is convenient to split things into
the three regions:
\begin{align}
        \mathcal{R}_1 \ &= \ \{r^* < -\frac{1}{2}s\} \ ,
    &\mathcal{R}_2 \ &= \ \{|r^*|\leqslant \frac{1}{2}s\}
    &\mathcal{R}_3 \ &= \ \{\frac{1}{2}s \leqslant r^*\} \ . \notag
\end{align}\ret

In the region $\mathcal{R}_1$ the bound on the
\eqref{end_linfty_bound}
is completely trivial. Here
the non-linear interaction breaks down entirely. All one has to do is
to use the bound:
\begin{align}
        (1-\frac{2M}{r}) \ &\lesssim  \ e^{\frac{1}{2M}r^*} \ ,
        &r^*\ \leqslant \ 0 \ , \notag
\end{align}
and the fact that the $\phi$ are uniformly bounded (which is the best
we can do!).\\

In the transition region $\mathcal{R}_2$,
we use the fine bound contained on the right
hand side of \eqref{global_sob}. This gives us that:
\begin{equation}
        (1-\frac{2M}{r})^\frac{1}{2}|\phi|^2 \ \lesssim \
        (1 + s)^{-2} \ . \notag
\end{equation}
Taking on the weight $(1 + |\bu| + |u|)$ which is $\sim(1+s)$ in
this region, and using the fact that $|\phi|^{p-2}\lesssim
(1+s)^{-\frac{1}{2}(p-2)}$ in $\mathcal{R}_2$ we again have
\eqref{end_linfty_bound}. \\

Finally, in the Minkowski like region $\mathcal{R}_3$ we have from
\eqref{global_sob} the uniform decay estimate $|\phi|\lesssim (1 + t
+ r^*)^{-1}$. This easily implies \eqref{end_linfty_bound}. This
completes our proof of the bootstrapping estimate
\eqref{need_to_bound} and hence our demonstration of Theorem
\ref{NLW_thm}.
\end{proof}

\ret

\ret

\ret\ret


\end{document}